\documentclass[12pt]{amsart}
\usepackage{amsmath,amsthm,amsfonts,latexsym,epsfig}

\newtheorem{thm}{Theorem}[section]

\newtheorem{pro}[thm]{Proposition}

\newcommand{\Aut}{{\mathrm{Aut}}\;}
\newcommand{\mre}{\iota}
\newcommand{\rank}{\mathsf{rank}}
\newcommand{\al}{\alpha}
\newcommand{\be}{\beta}
\newcommand{\ga}{\gamma}
\newcommand{\si}{\sigma}

\newcommand{\ta}{\theta}
\newcommand{\de}{{\delta}}
\newcommand{\la}{{\lambda}}

\newcommand{\vep}{\varepsilon}

\newcommand{\ld}{\ldots}
\newcommand{\cd}{\cdots}
\newcommand{\bA}{\mathbf{A}}
\newcommand{\bB}{\mathbf{B}}
\newcommand{\bM}{\mathbf{M}}
\newcommand{\De}{\Delta}
\newcommand{\cC}{\mathcal{C}}
\newcommand{\cP}{\mathcal{P}}

\newcommand{\cY}{\mathcal{Y}}
\newcommand{\cB}{\mathcal{B}}
\newcommand{\cT}{\mathcal{T}}
\newcommand{\cV}{\mathcal{V}}
\newcommand{\cS}{\mathcal{S}}
\newcommand{\fS}{\mathfrak{S}}
\newcommand{\bbQ}{\mathbb{Q}}
\newcommand{\bbZ}{\mathbb{Z}}
\newcommand{\bbP}{\mathbb{P}}
\newcommand{\bbC}{\mathbb{C}}
\newcommand{\pa}{\partial}

\newcommand{\bzero}{\mathbf{0}}

\newcommand{\dmrjdel}[1]{}

\title{The KP hierarchy, branched covers, and triangulations}

\author[I.~P.~Goulden]{I.P. Goulden$^*$}
\author[D.~M.~Jackson]{D.M.~Jackson$^\dagger$}

\thanks{
${\hspace{-1ex}}^*$Department of Combinatorics and Optimization,
                                                University of Waterloo, Waterloo, Ontario, Canada.;    \\
${\hspace{.35cm}}$ \texttt{ipgoulden@uwaterloo.ca}}

\thanks{
${\hspace{-1ex}}^\dagger$Department of Combinatorics and Optimization,
                                                University of Waterloo, Waterloo, Ontario, Canada.;  \\
${\hspace{.35cm}}$ \texttt{dmjackso@math.uwaterloo.ca}}

\date{March 27, 2008}

\begin{document}
\maketitle

\begin{abstract}
The KP hierarchy is a completely integrable system
of quadratic, partial differential equations that generalizes the KdV hierarchy.
A linear combination of Schur
functions is a solution to the KP hierarchy if and only if its coefficients
satisfy the Pl\"ucker relations from geometry. We give a solution
to the Pl\"ucker relations involving products of variables marking
contents for a partition, and thus give
a new proof of a content product solution to the
KP hierarchy, previously given by Orlov and Shcherbin. 
In our main result, we specialize this content product solution to prove that
the generating series for a general class of transitive ordered factorizations
in the symmetric group satisfies the KP hierarchy. 
These factorizations appear in geometry as encodings of branched covers, and thus
by specializing our transitive factorization result, we are able to prove
that the generating series for two classes of branched covers satisfies the KP hierarchy.
For the first of these, the double Hurwitz series, this result has been previously given
by Okounkov. The second of these, that we call the $m$-hypermap series,
contains the double Hurwitz series polynomially, as the leading coefficient in $m$.
The $m$-hypermap series also specializes further, first to the series for hypermaps and  then
to the series for maps, both in an orientable surface. For the latter series, we apply one
of the KP equations to obtain a new and remarkably simple recurrence for triangulations
in a surface of given genus, with a given number of faces. This recurrence leads
to explicit asymptotics for the number of triangulations with given genus and number
of faces, in recent work by Bender, Gao and Richmond. 


\end{abstract}


\section{Introduction and background}\label{intro}

\subsection{The KP hierarchy}

The KP (Kadomtsev-Petriashvili) hierarchy is a completely integrable system
of quadratic, partial differential equations for an unknown
function $F$, that generalizes the KdV hierarchy. Over the last two
decades, there has been strong interest in the relationship between integrable
systems and moduli spaces of curves. This began with Witten's Conjecture~\cite{w} for
the KdV equations, proved by Kontsevich~\cite{ko} (and more recently by a number of
others, including~\cite{kl}). Pandharipande~\cite{p} conjectured that solutions to
the Toda equations arose in a related context, which was proved by Okounkov~\cite{o}, who
also proved the more general result that a generating series for what he
called \emph{double Hurwitz numbers} satisfies the KP-hierarchy.  Kazarian
and Lando~\cite{kl}, in their recent proof of Witten's conjecture, showed that
it is implied by Okounkov's result for double Hurwitz numbers. More
recently, Kazarian~\cite{ka} has given a number of interesting results
about the structure of solutions to the KP hierarchy.

To fix ideas, the first few of the countable equations in the KP hierarchy are
\begin{equation}\label{canonKP}
F_{2,2}-F_{3,1}+\frac{1}{12}F_{1^4}+\frac{1}{2}F_{1^2}^2 = 0, 		\\
\end{equation}
\begin{equation*}
F_{3,2}-F_{4,1}+\frac{1}{6}F_{2,1^3}+F_{1,1}F_{2,1} = 0,			\\
\end{equation*}
\begin{equation*}
F_{4,2}-F_{5,1}+\frac{1}{4}F_{3,1^3}-\frac{1}{120}F_{1^6}+F_{1^2}F_{3,1}
+\frac{1}{2}F_{2,1}^2-\frac{1}{8}F_{1^3}^2-\frac{1}{12}F_{1^2}F_{1^4} = 0,
\end{equation*}
where 
$$ F_{r^{a_r}, \ldots ,1^{a_1}} :=  \frac{\partial^{a_r + \cdots  + a_1}}
{\partial p_r^{a_r} \cdots \partial p_1^{a_1}} F.  $$

We have found~\cite{mjd} to be an accessible source for the KP hierarchy, and
note that the variables  $p_i$ that we use in this paper are related to the
variables used in~\cite{mjd} by $p_i=i\, x_i$, $i\ge 1$. Thus the partial derivatives
in our equations (with respect the $p_i$'s) are scaled differently from those
in~\cite{mjd} (which are with respect to the $x_i$'s).

There is a powerful characterization of the solutions to the
KP hierarchy that is particularly convenient from an
algebraic combinatorics or geometric point of view. This
is well-known in the integrable systems
literature (see, \textit{e.g.}, Ch.~10 of~\cite{mjd}), and concerns 
linear combinations of Schur symmetric functions (here the
variables $p_1,p_2,\ld$ are the power sum symmetric functions).
The characterization states that the coefficients in the linear combination of Schur functions
satisfy the Pl\"ucker relations from geometry. In
the remainder of Section 1 (Sections 1.2 -- 1.4), we
give the background required to state this characterization precisely (Theorem~\ref{KPsol}).

In Section~2, we apply the Schur function characterization to give an explicit solution to the KP
hierarchy (as Theorem~\ref{Phisol}) in terms  of products of
a countable set of variables indexed by the integers, recording \textit{content}.
Content is a combinatorial parameter for partitions that often appears in
expressions for symmetric functions (which are themselves naturally indexed by
partitions).  We establish this result as an immediate corollary of
Theorem~\ref{Plucksol}, in which we give a content product solution
to the Pl\"ucker relations. The content product solution to the KP hierarchy
has previously been given by Orlov and Shcherbin~\cite{os}, using different methods.
In Section~3 we turn to algebraic combinatorics. We consider
a general set of transitive ordered factorizations in the symmetric group,
and prove that a particular generating function for
the numbers of such factorizations satisfies the KP hierarchy. This is
our main result, given as Theorem~\ref{T:main}. In Section~4,
we consider the geometrical interpretation
of transitive ordered factorizations, in terms of branched covers.
Thus in Section~4.2, as a corollary of our main result, we obtain Okounkov's~\cite{o} result, that the
generating series for double Hurwitz numbers satisfies
the KP hierarchy (Theorem~\ref{doubHurKP}). In Section~4.3 we consider a new class of geometric
numbers, called $m$-hypermap numbers, and prove that the generating series for
these satisfies the KP hierarchy (Theorem~\ref{mhypKP}). In Section 4.4,
we prove that the $m$-hypermap numbers are polynomials in $m$, with
the double Hurwitz numbers as the leading coefficient (Theorem~\ref{polym}),
and speculate that the rich geometry for the latter
might extend to the former. Then in Section 5, we specialize further to establish
that generating series for the numbers of rooted hypermaps (Theorem~\ref{hypKP}) and
rooted maps in an orientable surface (Theorem~\ref{mapKP}) are solutions to the
KP hierarchy. In Section 5.3, we conclude by using one of the KP equations to give a
new and remarkably simple recurrence for triangulations in an orientable surface
of specified genus, with a given number of faces (Theorem~\ref{trirec}).


\subsection{The Pl\"ucker relations}
If  $\la_1,\ld ,\la_n$ are integers with $\la_1\ge \cd\ge\la_n\ge 1$ and $\la_1+\cd +\la_n=d$,  then $\la=(\la_1,\ld ,\la_n)$ is said to be a \emph{partition} of $|\la|:=d$ (indicated by writing $\la\vdash d$) with $l(\la):=n$ \emph{parts}. The empty list $\vep$ of integers is to be regarded as a partition of $d=0$ with $n=0$ parts, and let $\cP$ denote the set of all partitions. If $\la$ has $f_j$ parts equal to $j$ for $j=1,\ld ,d$, then we also write $\la=(d^{f_d},\ld ,1^{f_1})$, where convenient. Also, $\Aut\la$ denotes the set of permutations of the $n$ positions that fix $\la$; therefore  $|\Aut\la |=\prod_{j\ge 1}f_j!$.

We consider a set $\{b_{\la}\in\bbQ[u_1,u_2,\ld]:\la\in\cP\}$, where $(u_1,u_2,\ld )$ is a list of indeterminates  independent of $p_1,p_2,\ld$. It is convenient to adopt two conventions for evaluating $b_{\la}$ in certain cases when $\la$ is not a partition. 

\noindent \textbf{Convention~1:\quad} For $\la=(\la_1,\ld ,\la_n)$, where $\la_1\ge \cd\ge\la_n$ and $\la_n\le 0$, then
\begin{equation}\label{conv1}
b_{\la}:=0,\qquad\text{if $\la_n<0$};\qquad\qquad\qquad
b_{\la}:=a_{(\la_1\ld ,\la_{n-1})},\qquad\text{if $\la_n=0$}.  
\end{equation}

The second convention concerns arbitrary lists $\la=(\la_1,\ld ,\la_n)$ of
integers (\textit{i.e.}, not necessarily in weakly decreasing order). For
such a list, we define the operator $\De_j$, for each $j=1,\ld ,n-1$, by
\begin{equation}\label{Delta}
\De_j\,\la=(\la_1,\ld,\la_{j-1},\la_{j+1}-1,\la_j +1,\la_{j+2},\ld,\la_n).
\end{equation}

\noindent \textbf{Convention~2:\quad}
\begin{equation}\label{conv2}
b_{\la}:=-b_{\De_j\,\la},\qquad j=1,\ld ,n-1.
\end{equation}
Note that if $\De_j\,\la=\la$, which is equivalent
to $\la_{j+1}=\la_j+1$, then~(\ref{conv2}) implies
that
\begin{equation}\label{fixed}
b_{\la}=-b_{\la}=0.
\end{equation}
It is now straightforward to determine the value of $b_{\la}$ for
any list $\la=(\la_1,\ld ,\la_n)$ of integers: if $\la$ is not in
weakly decreasing order, then we apply~(\ref{conv2}) until we
either can apply~(\ref{fixed}), or until we have a list that is
weakly decreasing; then we finish with~(\ref{conv1}) to
remove any terminal $0$'s.

Let $\al,\be\in\cP$, where $\al$ and $\be$ have $i$ and $j$ parts,
respectively, for some $i,j\ge 0$, $(i,j)\neq (0,0)$,
and suppose $\al =(\al_1,\ld ,\al_i)$ and $\be=(\be_1,\ld ,\be_j)$.
Let $m=\max\{ i+1,j-1,2\}$, and set $\al_{i+1}=\cd =\al_{m-1}$
$=\be_{j+1}=\cd =\be_{m+1}=0$. Then we say
that $\{ b_{\la}\in\bbQ[u_1,u_2,\ld]:\la\in\cP\}$ satisfies
the \textit{Pl\"ucker relations} if it satisfies the equation
\begin{equation}\label{plucker}
\sum_{k=0}^m(-1)^kb_{(\al_1-1,\ld,\al_{m-1}-1,\be_{k+1}+m-k)}
\,\cdot\, b_{(\be_1+1,\ld,\be_k+1,\be_{k+2},\ld ,\be_{m+1})} =0,
\end{equation}
for each such pair of partitions $\al,\be$, subject to~(\ref{conv1})
and~(\ref{conv2}). 
Note that partitions are represented by \textit{Maya diagrams} in~\cite{mjd}, and the statement of the Pl\"ucker relations given above is a translation of the Maya diagram notation used in~\cite{mjd}.

\subsection{Schur functions and characters of the symmetric group}
We recall a number of basic facts about symmetric functions (see~\cite{s}), where $x_1,x_2,\ld$ are the underlying indeterminates.
The $i$th \emph{power sum} symmetric function is  $p_i=\sum_{j\ge 1}x_j^i$, for $i\ge 1$, with  $p_0:=1$. The  $i$th \emph{complete} symmetric function $h_i$ is defined by 
$\sum_{i\ge 0}h_it^i=\prod_{j\ge 1}(1-x_jt)^{-1}$ and is related to the power sums through
\begin{equation}\label{hwithp}
\sum_{i\ge 0}h_it^i=\exp \sum_{k\ge 1}\frac{p_k}{k} t^k.
\end{equation}
The Schur function $s_{\la}$ 
may be expressed in terms of the complete symmetric functions  through the \emph{Jacobi-Trudi formula}
\begin{equation}\label{jactru}
s_{\la}=\det
\left( h_{\la_i-i+j}\right)_{i,j=1,\ld ,n},
\end{equation}
(with the convention that $h_i=0$ for $i<0$).
We write $h_i(p_1,p_2,\ld )$ to denote the $i$th complete symmetric function of indeterminates for
which the power sums are given by $p_1,p_2,\ld $, and similarly for  $s_{\la}(p_1,p_2,\ldots )$.

%
Let $\chi^{\la}_{\mu}$ denote the \emph{character} of the irreducible representation of $\fS_d$ indexed by $\la$, evaluated at any element of the  conjugacy class $\cC_{\mu}$ (we usually refer to $\chi^{\la}_{\mu}$ informally as
an irreducible character). Then, for $\la\vdash d$, the explicit expressions expressing the Shur functions and the power sums in terms of each other are
\begin{equation}\label{psandsp}
p_{\la}=\sum_{\mu\vdash d}\chi^{\mu}_{\la}s_{\mu},
\qquad\qquad   s_{\la}=\sum_{\mu\vdash d} \frac{|\cC_{\mu}|}{d!}\chi^{\la}_{\mu}p_{\mu}.
\end{equation}
It is convenient to consider a particular scaling of the irreducible characters, given by
\begin{equation}\label{scalechi}
g^{\la}_{\mu}=\vert\cC_{\mu}\vert\frac{\chi^{\la}_{\mu}}{\chi^{\la}_{(1^d)}},
\qquad\qquad \la ,\mu\vdash d.
\end{equation}
The following enumerative result
is well-known, and is included here since it will be applied later.

\begin{pro}\label{multperm}
For $\mu_i\vdash d\ge 1$ and $\si_i\in\cC_{\mu_i}$, $i=1,\ld ,k$, the
number of $k$-tuples $(\si_1,\ld ,\si_k)$, that satisfy the
equation $\si_1\cd\si_k=\mre$ (where $\mre$ is the identity permutation) is
$$\frac{1}{d!}\;\sum_{\la\vdash d}\left(\chi^{\la}_{(1^d)}\right)^2 g^{\la}_{\mu_1}\cd g^{\la}_{\mu_k}.$$
\end{pro}

\subsection{A characterization of solutions to the KP hierarchy}
For any solution $F$ to the KP hierarchy,
the series $e^F$ is called a $\tau$-\emph{function} of the KP hierarchy. The following result gives a
characterization of $\tau$-functions that is well-known in the integrable systems
literature (see, \textit{e.g.}, Ch.~10 of~\cite{mjd}).
\begin{thm}\label{tauschur}
The series $$\sum_{\la}b_{\la}s_{\la}(p_1,p_2,\ld )$$ is a $\tau$-function for the KP hierarchy if and only if $\{b_{\la}\in\bbQ[u_1,u_2,\ld]:\la\in\cP\}$ satisfies the Pl\"ucker relations~(\ref{plucker}).
\end{thm}
An obviously equivalent statement is the following.
\begin{thm}\label{KPsol}
The series $$\log\sum_{\la}b_{\la}s_{\la}(p_1,p_2,\ld )$$ is a solution to the KP hierarchy if and only if $\{b_{\la}\in\bbQ[u_1,u_2,\ld]:\la\in\cP\}$ satisfies the Pl\"ucker relations~(\ref{plucker}).
\end{thm}
The coefficients $b_{\la}$ in Theorems~\ref{tauschur} and~\ref{KPsol} are
called the \textit{Pl\"ucker coordinates} of the corresponding solution to the KP hierarchy.

\section{A content product solution to the KP hierarchy}

\subsection{Content products for partitions}
Some preliminary results about Ferrers diagrams of partitions are required.
The \emph{Ferrers graph} for a partition $\la$ is an array of unit boxes, called its \emph{cells}, with the $i$th row from the top containing $\la_i$ boxes, in columns (indexed from the left) $1,\ld ,\la_i$, for $i=1,\ld ,n$.  The \emph{content}  $c(w)$ of the cell $w$ in row $i$ and column $j$ is $c(w ):=j-i$.  For indeterminates $y_i$, where $i$ is an arbitrary integer, the \emph{content product} for $\la$ is
\begin{equation}\label{contprod}
C(\la):=\prod_{w\in\la}y_{c(w )},
\end{equation}
where the the product is over all cells  $w$ of the Ferrers graph of $\la$.
For example,  $C(5,3,3,2)=y_{-3}y_{-2}^2y_{-1}^2y_0^3y_1^2y_2y_3y_4$.
It is to be remembered throughout that the indeterminates $y_i$ may have negative suffices.
Perhaps the best known formula involving contents is
\begin{equation}\label{princspec}
\left. s_{\la}\right|_{p_i=x,i\ge 1} =\frac{d!}{\chi^{\la}_{(1^d)}}\prod_{w\in\la}(x+c(w)),
\end{equation}
known as the \textit{principal specialization} of the Schur
function (see, \textit{e.g.},~\cite{s}), which has been
recorded here since it will be applied later.


In addition, it will be convenient to consider other products of the $y_i$'s.
For pairs of integers $m,k$, we define $Y(m,k)$ by
\begin{equation}\label{Yprod}
Y(m,k):=
\begin{cases}
\prod_{j=1}^ky_{m+1-j},\qquad &\text{if $k\ge 1$,}\\
1,\qquad &\text{if $k=0$,}\\
Y(m-k,-k)^{-1},\qquad &\text{if $k\le -1$.}
\end{cases}
\end{equation}
Clearly, for all integers $m,j,k$ we have 
\begin{equation}\label{canc}
\frac{Y(m,k)}{Y(m,j)}=Y(m-j,k-j)=\frac{1}{Y(m-k,j-k)}.
\end{equation}


\begin{pro}\label{invar}
Let 
\begin{equation}\label{cY}
\cY (\la):=\prod_{i=1}^nY(\la_i-i+1,\la_i).
\end{equation}
where $\la$ is a list of integers of length $n\ge 2$. Then,
for $\De_j$ defined in~(\ref{Delta}), we have
$$\cY(\De_j\,\la)=\cY(\la),\qquad j=1,\ld ,n-1.$$
\end{pro}

\begin{proof}
{}From~(\ref{cY}),~(\ref{Delta}) and~(\ref{canc}), we have
\begin{eqnarray*}
\frac{\cY(\De_j\,\la)}{\cY(\la)}&=&
\frac{Y(\la_{j+1}-j,\la_{j+1}-1)Y(\la_j-j+1,\la_j+1)}{Y(\la_{j}-j+1,\la_{j})Y(\la_{j+1}-j,\la_{j+1})}\\
&=&\frac{Y(\la_j-j+1,\la_j+1)}{ Y(\la_j-j+1,\la_j)}
\frac{Y(\la_{j+1}-j,\la_{j+1}-1)}{Y(\la_{j+1}-j,\la_{j+1})}\\
&=&Y(-j+1,1)\frac{1}{Y(-j+1,1)} = 1,
\end{eqnarray*}
giving the result.
\end{proof}

Note that if $\la$ is a partition then 
\begin{equation}\label{contY}
\cY (\la )=C(\la).
\end{equation}

\subsection{A content product solution to the Pl\"ucker relations}

The following result gives an explicit class of solutions to the Pl\"ucker relations,
involving the content product $C(\la)$ defined in~(\ref{contprod}). We have been unable to find this
result stated explicitly in the literature.

\begin{thm}\label{Plucksol}
$$\{s_{\la}(q_1,q_2,\ld )\prod_{w\in\la}y_{c(w )}\, :\,\la\in\cP\}$$
satisfies the Pl\"ucker relations.
\end{thm}

\begin{proof}
For an arbitrary list $\la$ of integers, we define $f_{\la}$ by $f_{\la}:=\cY(\la)\,s_{\la}(q_1,q_2,\ld )$, through~(\ref{jactru}) and~(\ref{cY}).
We prove that $\{ f_{\la}:\la\in\cP\}$ satisfies
the Pl\"ucker relations~(\ref{plucker}). Here the indeterminates for
the $f_{\la}$'s are $y_i$ for integers $i$,
together with the $q_j$'s for positive integers $j$.
Note that the $f_{\la}$'s satisfies conventions~(\ref{conv1}), immediately
from~(\ref{jactru}), and~(\ref{conv2}) immediately from Proposition~\ref{invar}
and~(\ref{jactru}).


Local to this proof, we introduce the notation
\begin{equation*}
\ga^{(k)}:=(\al_1-1,\ld,\al_{m-1}-1,\be_{k+1}+m-k),
\qquad \de^{(k)}:=(\be_1+1,\ld,\be_k+1,\be_{k+2},\ld ,\be_{m+1}),
\end{equation*}
for $k=0,\ld ,m$,
and 
\begin{equation*}
\al':=(\al_1,\ld ,\al_{m-1}),\qquad\mbox{and}\qquad \be':=(\be_1,\ld ,\be_{m+1}).
\end{equation*}
Then, for $k=0,\ld ,m$, we obtain, applying~(\ref{cY}) and~(\ref{canc}),
\begin{eqnarray*}
\cY(\ga_k)\cY(\de_k)&=&\cY(\al'-1)\cY(\be'+1)
\frac{\cY(\be_{k+1}-k+1,\be_{k+1}+m-k)}{\cY(\be_{k+1}-k+1,\be_{k+1}+1)}
\prod_{i=k+2}^{m+1}\frac{\cY(\be_i-i+2,\be_i)}{\cY(\be_i-i+2,\be_i+1)}\\
&=&\cY(\al'-1)\cY(\be'+1)
\frac{\cY(-k,m-k+1)}{\prod_{i=k+2}^{m+1}\cY(-i+2,1)}\\
&=&\cY(\al'-1)\cY(\be'+1)y_{1-m}^{-1},
\end{eqnarray*}
where $\al'-1$ is the list obtained from $\al'$ by subtracting $1$ from every entry, and $\be'+1$ is the list obtained from $\be'$ by adding $1$ to every entry. In particular, we have
proved that $\cY(\ga_k)\cY(\de_k)$ is independent of $k$, so in checking the Pl\"ucker
relations~(\ref{plucker}), we have
\begin{equation}\label{pluchk}
\sum_{k=0}^m(-1)^k f_{\ga^{(k)}}\,\cdot\, f_{\de^{(k)}}=
\cY(\al'-1)\cY(\be'+1)y_{1-m}^{-1}\sum_{k=0}^m(-1)^k s_{\ga^{(k)}}\,\cdot\, s_{\de^{(k)}}.
\end{equation}

Now, for $m\ge 1$ consider the matrices
$$\bA=\left( h_{\al_i-i+j-1}\right)_{\left\{\substack{i=1,\ld ,m-1\\j=1,\ld ,m}\right.},
\qquad\mbox{and}\qquad
 \bB=\left( h_{\be_i-i+j+1}\right)_{\left\{\substack{i=1,\ld ,m+1\\j=1,\ld ,m}\right.},$$
where the $h_k$ are in the power sums $q_1,q_2,\ld$. Let $\bM$ be the $2m\times 2m$ matrix
given in the following partitioned form:
$$\bM=
\left(
\begin{tabular}{c|c}
$\bA$&$\bzero$\\ \hline
$\bB$&$\bB$\\
\end{tabular}
\right) ,
$$
where $\bzero$ is an $(m-1)\times m$ zero matrix. Then $\rank(\bM)\le 2m-1$, so $\det (\bM)=0$, and using the Laplace expansion with the columns partitioned into $\{ 1,\ld ,m\}$ and $\{ m+1,\ld ,2m\}$,
we obtain
$$0=\det (\bM) =\sum_{k=0}^m(-1)^k s_{\ga^{(k)}}\,\cdot\, s_{\de^{(k)}},$$
from~(\ref{jactru}). Together
with~(\ref{pluchk}),~(\ref{contY}) and~(\ref{contprod}), this implies the result.
\end{proof}

\subsection{A solution to the KP hierarchy}
As an immediate corollary to the content product solution for
the Pl\"ucker relations given in Theorem~\ref{Plucksol}, we now
give  a content product solution to the KP hierarchy. This result has
been previously obtained using different methods by Orlov and
Shcherbin~\cite{os} (see also Orlov~\cite{or}, equation~(1.19)).
\begin{thm}\label{Phisol}
The series
\begin{equation}\label{Phidef}
\Phi:=\log\sum_{\la}\left(\prod_{w\in\la}y_{c(w )}\right)s_{\la}(q_1,q_2,\ld )s_{\la}(p_1,p_2,\ld ).
\end{equation}
is a solution to the KP hierarchy (in the variables $p_1,p_2,\ld$).
\end{thm}

\begin{proof}
The result is immediate from Theorem~\ref{KPsol} and Theorem~\ref{Plucksol}.
\end{proof}

\section{Transitive ordered factorizations and the main result}

In this section we consider the
following general set of transitive ordered factorizations of permutations.
For $a_1, a_2, \ldots \ge 0$ and partitions $\al$ and $\be$ of $d\ge 1$, let $\cB^{(a_1,a_2, \ldots)}_{\al ,\be}$ be the set of tuples of
permutations $(\si,\ga, \pi_1,\pi_2,\ld )$ on $\{ 1,\ld ,d\}$ such that

\begin{description}
\item[C1] $\si\in\cC_{\al},$ $ \ga\in\cC_{\be}$,  and $d-l(\pi_i) = a_i$ for $i\ge1$, where $l(\pi_i) $ is the number of cycles in the disjoint cycle decomposition of $\pi_i$;
\item[C2] $\si\ga\pi_1\pi_2\cdots =\mre$;
\item[C3] $\langle \si,\ga, \pi_1,\pi_2,\ld\rangle$ acts transitively on $\{ 1,\ld ,d\}$.
\end{description}
Let $b^{(a_1,a_2, \ldots)}_{\al ,\be}$ be the number of tuples in $\cB^{(a_1,a_2, \ldots)}_{\al ,\be}$, and
let ${\widetilde{b}}^{(a_1,a_2, \ldots)}_{\al ,\be}$ be the number of these tuples in the case in which  condition~{\bf C3} is not invoked.
Instances of such transitive factorizations appear in the combinatorial literature
in many places (see, \textit{e.g.},~\cite{lz}, where they are called \textit{constellations}).

As a corollary of Theorem~\ref{Phisol},
we now prove our main result, that a particular generating series for the
numbers $b^{(a_1,a_2, \ldots)}_{\al ,\be}$ of transitive ordered factorizations
is a solution to the KP hierarchy.

\begin{thm}\label{T:main}
The series
\begin{equation}\label{Bbseries}
B :=  \sum_{\substack{|\al|=|\be|=d\ge1,\\ a_1,a_2,\cd \ge0}}\frac{1}{d!} \,
b^{(a_1,a_2,\ldots)}_{\al,\be} p_\al q_\be \, u_1^{a_1} u_2^{a_2} \cdots
\end{equation}
is
a solution to the KP hierarchy (in the variables $p_1,p_2,\ld).$
\end{thm}
\begin{proof}
From the exponential formula for exponential generating
series (see, \textit{e.g.},~\cite{gj1}; interestingly, Hurwitz seems to have been the
first person to write this down clearly), we have
\begin{equation}\label{e:CtC}
B= \log(\tilde B)
\end{equation}
where
$$
\tilde{B} := 1+
\sum_{\substack{|\al|=|\be|=d\ge1,\\ a_1,a_2,\cd \ge0}}\frac{1}{d!}
\tilde{b}^{(a_1,a_2,\ldots)}_{\al,\be} p_\al q_\be \, u_1^{a_1} u_2^{a_2} \ldots.
$$
But, from Proposition~\ref{multperm}, (\ref{psandsp}) and~(\ref{scalechi}) we obtain
\begin{align*}
\tilde B 
&=
\sum_{\al,\be,\la\vdash d\ge0}  \frac{1}{d!^2}
\left(\chi^\la_{(1^d)}\right)^2 p_\al\,q_\be\,
g^\la_\al\, g^\la_\be
\prod_{i \ge 1} \left(\sum_{\mu_i\vdash d}
g^\la_{\mu_i}\, u_i^{d-l(\mu_i)}\right) \\
&=\sum_{\la\vdash d\ge0} s_\la(q_1,q_2,\ldots)\, s_\la (p_1, p_2,\ld)
\prod_{i \ge 1} \left(\sum_{\mu_i\vdash d}
g^\la_{\mu_i}\, u_i^{d-l(\mu_i)}\right).
\end{align*}
It is now immediate from~(\ref{psandsp}) and~(\ref{princspec}) that
\begin{align*}
\sum_{\mu_i\vdash d} g^\la_{\mu_i}\, u_i^{d-l(\mu_i)}
= \left. \frac{d! u_i^d}{\chi^\la_{(1^d)}} s_\la  \right|_{p_j=u_i^{-1}, \, j\ge1} 
= u_i^d \, \prod_{w\in\la} \left(u_i^{-1} + c(w)\right)
=\prod_{w\in\la} (1+u_i\, c(w) ).
\end{align*}
This, together with~(\ref{Phidef}),~(\ref{Bbseries}) and~(\ref{e:CtC}), gives
\begin{equation}\label{PhiB}
B = \left. \Phi \right|_{y_j=\prod_{i\ge1}(1+j\,u_i),\, j\in\bbZ }
\end{equation}
and the result now follows from Theorem~\ref{Phisol}. Note that
\begin{equation}\label{elemk}
\prod_{i\ge1}(1+j\,u_i) = 1+ \sum_{k\ge 1} e_k(u_1, u_2, \ld)\,j^k,
\end{equation}
where $e_k(u_1, u_2, \ld)$ is the $k$th \textit{elementary} symmetric function
in $u_1, u_2, \ld$.
\end{proof}


\section{Branched covers, double Hurwitz numbers and $m$-hypermap numbers}

\subsection{Branched covers}
The transitive ordered factorizations in  $\cB^{(a_1,a_2, \ldots)}_{\al ,\be}$ also
have geometric significance for, by an encoding due
to Hurwitz~\cite{hur}, they  correspond to branched covers of $\bbC\,\bbP^1$ with fixed branched
points, say $0$, $\infty$, and $X_i$, $i\ge 1$. We require the branching
over $0$ and $\infty $ to be $\al$ and $\be$, respectively, and
the branching over $X_i$ to have $d-a_i$ cycles in the disjoint cycle decomposition, $i\ge 1$.
Thus, in {\bf [C1]}, the permutations  $\si$ and $\ga$ encode the branching over $0$ and $\infty$,
respectively, and $\pi_i$ encodes the branching over $X_i$, $i\ge 1$.
{\bf [C2]} is a monodromy condition
and {\bf [C3]} makes the covers connected.
The genus $g$ of these branched covers follows from the Riemann-Hurwitz formula,
which in this case gives
\begin{equation}\label{RHg}
a_1+a_2+\cdots=r^g_{\al,\be},
\end{equation}
where, for partitions $\al,\be$, and non-negative integer $g$,
\begin{equation}\label{rgalbe}
r_{\al,\be}^g=l(\al)+l(\be)+2g-2.
\end{equation}

\subsection{Double Hurwitz numbers}

Double Hurwitz numbers arise in the enumeration of branched covers (see~\cite{o},~\cite{gjv}),
where they correspond to transitive ordered factorizations in $\cB^{(a_1,a_2,\ldots )}_{\al ,\be}$,
in which the branching over each point $X_i$ is simple (a transposition) for $i=1,\ldots ,r^g_{\al ,\be}$,
and trivial (the identity permutation) for $i> r^g_{\al ,\be}$. Thus,
rescaled for geometric reasons, as in~\cite{gjv}, the
\textit{double Hurwitz number} $H^g_{\al ,\be}$ is defined by
\begin{equation}\label{Hgdef}
H^g_{\al ,\be}=\frac{1}{d!}\vert\Aut\,\al\vert\vert\Aut\,\be\vert
b^{(a_1,a_2,\ldots )}_{\al ,\be},
\qquad\qquad \al ,\be\vdash d\ge 1,\qquad g\ge 0,
\end{equation}
where $a_i=1$ for $i=1,\ldots ,r^g_{\al ,\be}$, and $a_i=0$ for $i> r^g_{\al ,\be}$,
and $g$ is defined by~(\ref{rgalbe}). Among the results known for double Hurwitz numbers,
there is the beautiful and explicit formula for the
case $\be=(1^d)$  and $g=0$ (see, \textit{e.g.},~\cite{gj2})
\begin{equation}\label{gjform}
H^{0}_{\al ,(1^d)}=d!\,
d^{\,l(\al )-3} (d+l(\al )-2)!
\prod_{i=1}^{l(\al)}\frac{\al_i^{\al_i}}{\al_i!}.
\end{equation}

We now prove that a particular generating series for double Hurwitz
numbers is a solution to the KP hierarchy, as a corollary of Theorem~\ref{T:main}.
This result was first proved by Okounkov~\cite{o}, using a different method,
and then more recently by Orlov~\cite{or}, and Kazarian~\cite{ka}.
\begin{thm}\label{doubHurKP}
The double Hurwitz series
\begin{equation}\label{douHseries}
H=\sum_{\substack{\vert\al\vert=\vert\be\vert\ge 1,\\ g\ge 0}}
\frac{H^g_{\al ,\be}}{\vert\Aut\,\al\vert\vert\Aut\,\be\vert}
p_{\al}q_{\be}\frac{t^{r_{\al,\be}^g}}{r_{\al,\be}^g!}
\end{equation}
is a solution to
the KP hierarchy (in the variables $p_1,p_2,\ld$).
\end{thm}

\begin{proof}
{}From~(\ref{Bbseries}),~(\ref{Hgdef}) and~(\ref{douHseries}), we obtain
$$\left[\frac{t^{r_{\al,\be}^g}}{r_{\al,\be}^g!}\right]H
= \left[u_1\cdots u_{r_{\al,\be}^g}\right] B.$$
But, from~\cite{gj1} (Lemma~4.2.5(1),  \textit{p.} 233), this implies that
$$H=\left. B\right|_{e_k(u_1,u_2,\ldots)=\frac{t^k}{k!},k\ge 1},$$
and the result follows from Theorem~\ref{T:main}. Note that, in terms of
the series $\Phi$, we have
\begin{equation}\label{HPhi}
H=\left.\Phi\right|_{y_j=e^{jt},j\in\bbZ},
\end{equation}
from~(\ref{PhiB}) and~(\ref{elemk}).
\end{proof}

\subsection{$m$-Hypermap numbers}

Let $m$ be a fixed positive integer.
Define $c^{(g,m)}_{\al ,\be}$ by
\begin{equation}\label{mhypsum}
c^{(g,m)}_{\al ,\be}:=\sum b^{(a_1,a_2,\ldots)}_{\al ,\be},
\end{equation}
where the sum is over all $(a_1,a_2,\ldots)$ with $a_i=0$ for $i>m$, and 
\begin{equation}\label{sumrestr}
a_1+\cdots +a_m=r_{\al,\be}^g.
\end{equation}
Thus we are considering genus $g$ branched covers with branching over $0$ and $\infty$ specified
by $\al$ and $\be$, respectively, and arbitrary branching at $m$ other points $X_1,\ldots ,X_m$.
For geometric reasons, we scale these numbers in the same way as for
double Hurwitz numbers above. Hence we
define the $m$-\textit{hypermap number} $N^{(g,m)}_{\al ,\be}$ by
\begin{equation}\label{mhypnum}
N^{(g,m)}_{\al ,\be}:=\frac{1}{d!}\vert\Aut\,\al\vert\vert\Aut\,\be\vert
c^{(g,m)}_{\al ,\be},
\qquad\qquad \al ,\be\vdash d\ge 1,\qquad g\ge 0.
\end{equation}
We use the term $m$-hypermap number because the case $m=1$ yields \textit{rooted hypermaps},
as discussed in a later section.
The case $\be=(1^d)$ and $g=0$ has been considered
by Bousquet-M\'elou and Schaeffer~\cite{bms}, where they obtained
the beautiful and explicit formula
\begin{equation}\label{bmsform}
N^{(0,m)}_{\al ,(1^d)}=d!\,
m\frac{\left((m-1)d-1\right)!}{\left( (m-1)d-l(\al )+2\right) !}
\prod_{i=1}^{l(\al)}\binom{m\al_i-1}{\al_i}.
\end{equation}

As a second corollary of Theorem~\ref{T:main}, we 
now prove that a particular generating series for $m$-hypermap
numbers is a solution to the KP hierarchy.
\begin{thm}\label{mhypKP}
The $m$-hypermap series
\begin{equation}\label{mhypseries}
N^{(m)}=\sum_{\substack{\vert\al\vert=\vert\be\vert =d\ge 1,\\ g\ge 0}}
\frac{N^{(g,m)}_{\al ,\be}}{\vert\Aut\,\al\vert\vert\Aut\,\be\vert}
p_{\al}q_{\be}t^{r_{\al,\be}^g}.
\end{equation}
is a solution to
the KP hierarchy (in the variables $p_1,p_2,\ld$).
\end{thm}

\begin{proof}
{}From~(\ref{Bbseries}),~(\ref{mhypnum}) and~(\ref{mhypseries}),  we obtain
$$N^{(m)}=\left. B\right|_{u_1=\cdots =u_m=t, u_j=0, j> m},$$
and the result follows from Theorem~\ref{T:main}. Note that, in terms of
the series $\Phi$, we have
\begin{equation}\label{NmPhi}
N^{(m)}=\left.\Phi\right|_{y_j=(1+jt)^m,j\in\bbZ},
\end{equation}
from~(\ref{PhiB}).
\end{proof}

\subsection{A direct relationship between Hurwitz numbers and $m$-hypermap numbers.}
One relationship between  Hurwitz numbers and $m$-hypermap numbers arises from
inclusion-exclusion, as follows. It is straightforward that the summation
on the right hand side of~(\ref{mhypsum}), subject to~(\ref{sumrestr}), and the
further restriction that none of $a_1,\ld ,a_m$ is equal to $0$, is given by
\begin{equation}\label{incexcm}
\sum_{j\ge 0}(-1)^j c^{(g,m-j)}_{\al ,\be}.
\end{equation}
But this forces $a_i\ge 1$ for all $i=1,\ld ,m$, and so in the case $m=r^g_{\al ,\be}$, we
must have $a_i=1$ for all $i=1,\ld ,m$. This implies that branching over each such $X_i$ is
simple, and so, rescaling as in~(\ref{Hgdef}) and~(\ref{mhypnum}), we have
$$\sum_{j\ge 0}(-1)^j N^{(g,r^g_{\al ,\be}-j)}_{\al ,\be}=H^g_{\al ,\be}.$$
This inclusion-exclusion argument was given 
in~\cite{bms} for the case $g=0$ and $\be=(1^d)$, and enabled them to obtain~(\ref{gjform})
from~(\ref{bmsform}).

The following result gives another relationship, perhaps more direct, obtained by
comparing the generating series $H$ and $N^{(m)}$ in~(\ref{douHseries}) and~(\ref{mhypseries}). 
\begin{thm}\label{polym}
For $g\ge 0$ and partitions $\al$ and $\be$ of $d\ge 1$, $N^{(g,m)}_{\al,\be }$ is
a polynomial in $m$ of degree $r^g_{\al ,\be }$, over $\bbQ$. Moreover,
$$\left[\frac{m^{r^g_{\al ,\be }}}{r^g_{\al ,\be }!}\right] N^{(g,m)}_{\al,\be }
=H^{g}_{\al,\be }.$$
\end{thm}

\begin{proof}
The result follows immediately by comparing~(\ref{HPhi}) and~(\ref{NmPhi}).
\end{proof}

For example, it is straightforward to check that Theorem~\ref{polym} holds
in the case $\be = (1^d)$ and $g=0$, using the explicit expressions
given for these particular double Hurwitz and $m$-hypermap numbers in~(\ref{gjform}) and~(\ref{bmsform}),
respectively.

We do not know an elementary direct proof of Theorem~\ref{polym}.
There is a remarkably rich literature on the geometry associated with Hurwitz numbers (the case $\be=(1^d)$ of
double Hurwitz numbers). The fact that Hurwitz numbers arise in Theorem~\ref{polym} as the leading coefficient of $m$-hypermap
numbers (where we can specialize to $\be=(1^d)$ in the same way) causes us to speculate that much of
the geometry associated with Hurwitz numbers may extend to $m$-hypermap numbers.


\section{Hypermaps, maps and triangulations in orientable surfaces}

\subsection{Hypermaps in orientable surfaces.}
A connected graph embedded in an orientable surface partitions the surface into regions called \emph{faces}, and  for \emph{two-cell} embeddings, which are considered here,  the faces are homeomorphic to open discs.  If the faces are properly \emph{two-colourable}, so faces of the same colour intersect only at vertices  (using colours black and white) the embedded graph is called a \emph{hypermap}, where the black faces are \emph{hyperedges} and the white faces are \emph{hyperfaces}. The \emph{degree} of a vertex is the number of adjacent hyperedges, and the degree of a hyperedge or hyperface is the number of sides of edges encountered when traversing the boundary once.  We consider \emph{rooted} hypermaps, in which  an arrow is drawn on one edge from a tail vertex to a head vertex, so that, moving around the tail vertex, there is a white face to the counterclockwise side of the root edge. 
\begin{figure}[ht]
\begin{center}
\scalebox{.60}{\includegraphics{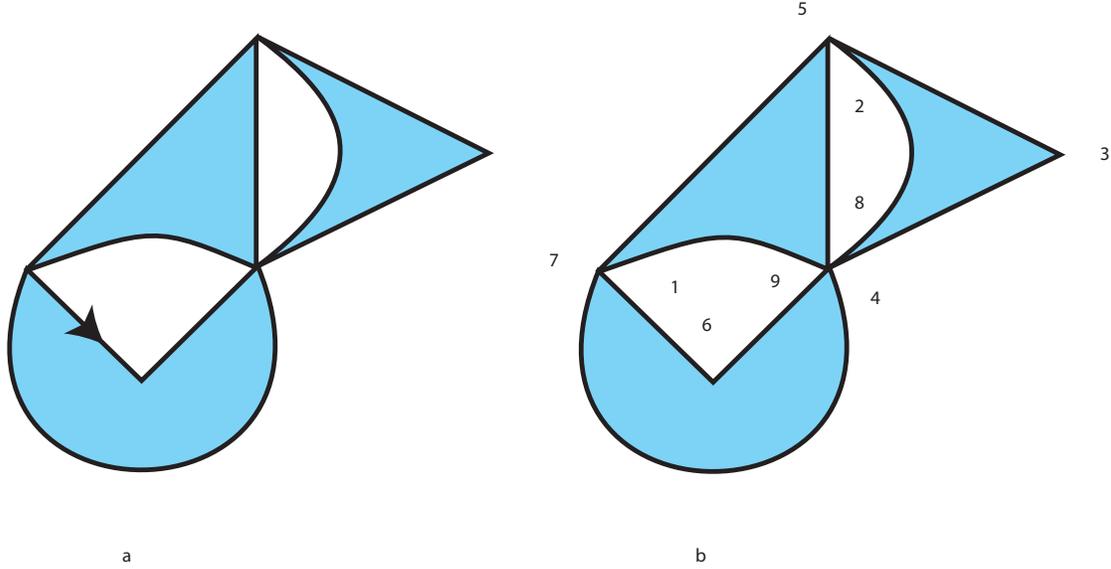}}
\end{center}
\caption{A rooted hypermap on the sphere.}\label{hypereg}
\end{figure}
For example, on the
left side of Figure~\ref{hypereg}, marked as ``a'', we give a rooted hypermap in the sphere, in which
there are three faces of each colour  (the black faces have shaded interiors).

Let $M^g_{\al ,\be}$ denote the number of rooted hypermaps in an orientable surface of genus $g$, with vertex-degrees specified by the parts of the partition $\al$, and hyperedge (black face) degrees specified by the parts of $\be$. It is well-known (see, \textit{e.g.},~\cite{jv1},~\cite{jv2},~\cite{lz},~\cite{t}) that 
\begin{equation}\label{hyprot}
M^g_{\al ,\be}=\frac{1}{(d-1)!}c^{(g,1)}_{\al ,\be},\qquad\qquad
\al,\be\vdash d\ge 1,g\ge 0.
\end{equation}
Indeed, this explains the use of the term $m$-hypermap numbers for $N^{(g,m)}_{\al ,\be}$ (which
is simply a rescaling of $c^{(g,m)}_{\al ,\be}$, \textit{via}~(\ref{mhypnum})).
The $(d-1)!:1$ correspondence implicit in~(\ref{hyprot}) is described as follows: Label the corners
of the white faces $1,\ld ,d$ so that $1$ is assigned to the corner on the
counterclockwise side of the root edge, as encountered when moving around the tail
vertex. The remaining $d-1$ labels may be placed arbitrarily, to give the
factor $(d-1)!$. From this labelled hypermap, we obtain three permutations
in $\fS_d$: $\si$, $\ga$ and $\pi$. Each disjoint cycle of $\si$ gives the
labels encountered when moving around a vertex in a counterclockwise direction.
Each disjoint cycle of $\pi$ gives the  labels encountered when traversing
the interior boundary of a hyperface in a counterclockwise  direction. Define
the label of a black corner (\textit{i.e.}, a corner of a hyperedge) to be
the label of the white corner that is encountered in the clockwise direction
when moving around their common vertex.
Based on this labelling convention, each disjoint cycle of $\ga$ gives the
labels encountered when traversing the interior boundary
of a hyperedge in a counterclockwise direction.
Given a rooted hypermap counted by $M^g_{\al ,\be}$, it is
clear by construction that $\si\in\cC_{\al}$ and
 $\ga\in\cC_{\be}$, and that $\si\ga\pi=\mre$. Moreover, $\langle \si ,\ga ,\pi\rangle$ acts
transitively on $\{ 1,\ldots ,d\}$, because the underlying graph is connected. This
is clearly reversible, and specifies the required $(d-1)!:1$ correspondence.

For example, on the right side of Figure~\ref{hypereg}, marked as ``b'', is one of the labellings of the
hypermap given on the left side. For this example, we obtain
$$\si=(1\, 7)(2\, 5)(3)(4\, 8\, 9)(6),\qquad
\ga=(1\, 8\, 5)(2\, 4\, 3)(6\, 7\, 9),\qquad
\pi=(1\, 6\, 9)(2\, 8)(3\, 4\, 7\, 5),$$
and it is easy to verify that, in this case, we have $\si\ga\pi =\mre$.

It is instructive to consider the genus $g$ in~(\ref{RHg}), which gives the genus
for the corresponding branched cover. In terms of the hypermap in the
present case (with $m=1$), condition~(\ref{RHg}) gives
\begin{equation}\label{coverg}
d-l(\pi)=l(\al)+l(\be)+2g-2.
\end{equation}
Now, in the terminology of Euler's polyhedral formula, the hypermap has $V=l(\al)$ vertices,  $E=d$ edges, and $F=l(\be)+l(\pi)$ faces, so a straightforward rearrangement of~(\ref{coverg}) gives
$$l(\al)-d+\left( l(\be)+l(\pi)\right)=2-2g,$$
or $V-E+F=2-2g$, which is Euler's formula. Our use of genus is therefore
consistent for covers and hypermaps. 

As a specialization of Theorem~\ref{mhypKP}, we now
prove that a particular generating series for the rooted hypermap
numbers $M^{g}_{\al ,\be}$ is
a solution to the KP hierarchy.

\begin{thm}\label{hypKP}
The rooted hypermap series
\begin{equation}\label{hmapser}
M:=\sum_{\substack{\vert\al\vert=\vert\be\vert=d\ge 1,\\ g\ge 0}}
\frac{M^{g}_{\al ,\be}}{d}
p_{\al}q_{\be}t^{r_{\al,\be}^g}
\end{equation}
is a solution to the KP hierarchy (in the variables $p_1,p_2,\ld $).
\end{thm}

\begin{proof}
{}From~(\ref{mhypnum}),~(\ref{mhypseries}) and~(\ref{hyprot}), we have $M=N^{(1)}$, and the result follows immediately, as the specialization of Theorem~\ref{mhypKP} to $m=1$.
\end{proof}

Note that the hypermap series $M$ enables us to record vertex degrees and hyperedge degrees separately  but not hyperface degrees, so only the total number of hyperfaces (through Euler's formula) is recorded. The analogous series in which hyperface degrees (in addition to vertex and hyperedge degrees) are recorded is \emph{not} a solution to the KP hierarchy.

\subsection{Maps in orientable surfaces} 
The specialization from rooted hypermaps to \emph{rooted maps} is by requiring each hyperedge to be of degree $2$, and then ``collapsing'' each hyperedge (black face) to a single edge where, for the
hyperedge containing the directed root edge, the collapsed single edge has the same direction as the directed edge. 
The result is a \textit{rooted map} with the collapsed singled edges as its edges and the hyperfaces as its faces. Let $R_{\al}^{(n,m)}$ denote the number of rooted maps in an orientable surface with $n$ edges, $m$ faces, and vertex degrees  specified by the parts of $\al$. Thus $\al$ is a partition of $2n$, and the genus $g$ of the surface, by Euler's formula, is given by $l(\al)-n+m=2-2g$. Note that, by duality, $R_{\al}^{(n,m)}$ is also equal to the number of rooted maps in an orientable surface with $n$ edges, $m$ vertices, and face degrees specified by the parts of $\al$.

As a specialization of Theorem~\ref{hypKP}, we now
prove that a particular generating series for the rooted map
numbers $R^{(n,m)}_{\al}$ is
a solution to the KP hierarchy.

\begin{thm}\label{mapKP}
The rooted map series 
\begin{equation}\label{mapser}
R:=\sum_{n,m\ge 1}\sum_{\al\vdash 2n}\frac{R_{\al}^{(n,m)}}{2n}p_{\al}w^mz^n
\end{equation}
is a solution to the KP hierarchy (in
the variables $p_1,p_2,\ld $).
\end{thm}

\begin{proof}
Comparing~(\ref{hmapser}) and~(\ref{mapser}) and
applying Euler's formula and~(\ref{rgalbe}), we
have $\be=(2^n)$ and $m=2n-r^g_{\al,\be}$, so
$$R=
\left. M\right|_{\left\{\substack{q_2=w^2z,t=w^{-1}\\q_i=0,i\neq 2}\right.},\; .$$
The result follows immediately from Theorem~\ref{hypKP}.
\end{proof}

For results that correspond to Theorems~\ref{hypKP} and~\ref{mapKP}
when the generating series are expressed as matrix models see, \textit{e.g.},~\cite{or}
(and for the connection between matrix models and generating series for
rooted maps in an orientable surface see, \textit{e.g.},~\cite{j}.) 

\subsection{Triangulations in an orientable surface of arbitrary genus}
In this final section, we apply Theorem~\ref{mapKP} to obtain a recurrence equation
for rooted \textit{cubic} maps (all vertices have degree $3$) in
an orientable surface. By duality, these
are equivalent to rooted \textit{triangulations}, an important class of maps for the study
of surfaces in general.

We begin by defining, for any $\mu\subseteq\{ 1,2,3\}$, the generating
series $V_{\mu}(p_{\mu},w,z)$ by
\begin{equation}\label{parRmu}
V_{\mu}(p_{\mu},w,z):=\left. 2z\frac{\pa}{\pa z}R\right|_{p_i=0,i\notin\mu}.
\end{equation}
In the next result, we obtain combinatorial relationships between various $V_{\mu}$.
\begin{pro}\label{cubeconst}
\begin{eqnarray*}
(i)&\;\;V_{\{1,2,3\}}(p_1,p_2,p_3,w,z)=\frac{p_2w^2z}{1-p_2z}+\frac{1}{1-p_2z}
V_{\{1,3\}}\left( p_1,p_3,w,\frac{z}{1-p_2z}\right) ,\\
(ii)&\;\;V_{\{1,3\}}(p_1,p_3,w,z)=wz^{-1}T^2+\frac{4p_1p_3 w^2z^2}{1-4p_1p_3z^2}+
\frac{1}{1-4p_1p_3z^2}V_{\{3\}}\left( p_3,w,\frac{z}{\sqrt{1-4p_1p_3z^2}}\right) ,
\end{eqnarray*}
where $T$ is the formal power series solution to
\begin{equation}\label{lagtree}
T=z(p_1+p_3T^2).
\end{equation}
\end{pro}

\begin{proof}
Comparing~(\ref{parRmu}) and~(\ref{mapser}), we obtain
$$V_{\mu}=\sum_{n,m\ge 1}\sum_{\al\vdash n} R_{\al}^{(n,m)}p_{\al}w^mz^n,$$
with the restriction in the sum over $\al$ that all parts of $\al$ are
contained in $\mu$. Thus, $V_{\mu}$ is the ordinary generating  series
for the set $\cV$ of rooted maps in an orientable surface in which all vertex degrees
are contained in $\mu$. 
\vspace{.05in}

\noindent\underline{\textit{For Part $(i)$}}: \quad There are two cases for the maps in $\cV_{\{ 1,2,3\}}\,$: 
\begin{description}
\item[Case 1] All vertices are of degree $2$; 
\item[Case 2] Some vertex has degree $1$ or $3$.
\end{description}

In Case 1, there is exactly one map with $k$ edges for each $k\ge 1$, namely
the  $k$-cycle embedded in the sphere.  This accounts for the first term
on the right hand side of Part $(i)$ of the result.

In Case 2, each such map in $\cV_{\{ 1,2,3\}}$ can be uniquely created
by subdividing the edges of maps in $\cV_{\{ 1,3\}}$, replacing them
by paths in which all internal vertices have degree $2$. The
number of faces is unchanged in this construction. This accounts for the second term
on the right hand side of Part $(i)$ of the result, where the external factor $(1-p_2z)^{-1}$ is
an adjustment for the root edge.
\vspace{.1in}

\noindent\underline{\textit{For Part $(ii)$}}:  \quad There are three cases for the maps in $\cV_{\{ 1,3\}}$: 
\begin{description}
\item[Case 1] The map has $0$ cycles;  
\item[Case 2] The map has $1$ cycle; 
\item[Case 3] The map has at least $2$ cycles. 
\end{description}

In Case 1, these maps are rooted trees in the sphere, in which all vertices have
degree $1$ or $3$ (each of these has $1$ face). If we ``\textit{cut}'' the root edge, then such trees
decompose into an ordered pair of trees from $\cT$, which consists of rooted
ordered trees on at least one vertex, in which every vertex has up-degree $0$ or $2$.
Let $T_{i_1,i_3}^{n}$ denote the number of trees in $\cT$ with $i_1$ vertices
of up-degree $0$, $i_3$ vertices of up-degree $2$, and $n$ edges (so $n=i_1+i_3-1$),
and let
$$T:=\sum_{\substack{n,i_1,i_3\ge 0\\ n=i_1+i_3-1}}
T_{i_1,i_3}^{n}p_1^{i_1}p_3^{i_3}z^{n+1}.$$
Then, $T$ clearly satisfies the functional equation~(\ref{lagtree}), and so
we obtain the first term
on the right hand side of Part $(ii)$ of the result.

In Case 2, these are embedded in the sphere, with $2$ faces, and each vertex
on the cycle has degree $3$. The edge incident
with this vertex that does not lie on the cycle is either inside or
outside the cycle, and is also incident with the root vertex of
a tree in $\cT$.  Thus the contribution to $V_{\{1,3\}}$ in this
case is
\begin{equation}\label{case2T}
w^2\frac{2p_3zT}{1-2p_3zT}+w^2 \frac{2p_3z\frac{z\pa}{\pa z}T}{1-2p_3zT},
\end{equation}
where the first term in~(\ref{case2T}) is for maps with root edge
on the cycle (in a canonical, say clockwise direction), and
the second term in~(\ref{case2T}) is for maps with root edge off
the cycle. For this second term, we place the root edge on some canonical
side of the cycle, say outside, and then direct it in either of
the $2$ possible directions.
But, applying $\frac{z\pa}{\pa z}$ to~(\ref{lagtree}), we obtain
\begin{equation}\label{zddzT}
\frac{z\pa}{\pa z}T=\frac{T}{1-2p_3zT},
\end{equation}
and solving~(\ref{lagtree}) as a quadratic equation in $T$, we obtain
$$T=\frac{1-\sqrt{1-4p_1p_3z^2}}{2p_3z},$$
where we have rejected the other root since it is not a formal power series. 
But this explicit expression for $T$ gives
\begin{equation}\label{denomT}
\frac{1}{1-2zp_3T}=\frac{1}{\sqrt{1-4p_1p_3z^2}}.
\end{equation}
Simplifying~(\ref{case2T}) by means of~(\ref{zddzT}) and~(\ref{denomT}), we
obtain the second term on the right hand side of Part $(ii)$ of the result.

In Case 3, each such map in $\cV_{\{ 1,3\}}$ can be uniquely created
by subdividing the edges of maps in $\cV_{\{ 3\}}$, replacing them
by paths in which each internal vertex has degree $3$. The edge incident with
this vertex that does not lie on the path is in the face on either side of the path,
and is also incident with the root vertex of a tree in $\cT$. The
number of faces is unchanged in this construction. The contribution to $V_{\{1,3\}}$ in this
case is
$$ \left(\frac{1}{1-2zp_3T}+\frac{2p_3z\frac{z\pa}{\pa z}T}{1-2zp_3T}\right)
V_{\{3\}}\left( p_3,w,\frac{z}{1-2zp_3T}\right)  ,$$
where the external factor is
an adjustment for the root edge. Simplifying by means
of~(\ref{zddzT}) and~(\ref{denomT}), we
obtain the third term
on the right hand side of Part $(ii)$ of the result.
\end{proof}

Now let $\cS=\{ (n,g)\in\bbZ\times\bbZ : n\ge -1, 0\le g\le\frac{n+1}{2}\}$, and
define $f(n,g)$ by the quadratic recurrence equation
\begin{equation}\label{trirecence}
f(n,g):=\frac{4(3n+2)}{n+1}
\left(n(3n-2)f(n-2,g-1)+\sum f(i,h)f(j,k)\right),
\end{equation}
for $(n,g)\in\cS\setminus\{(-1,0),(0,0)\}$, where the summation
is over $(i,h)\in\cS,(j,k)\in\cS$ with $i+j=n-2$ and $h+k=g$,
subject to the initial conditions
\begin{equation*}
f(-1,0)=\frac{1}{2},\qquad\qquad\qquad f(n,g)=0, \qquad (n,g)\notin\cS.
\end{equation*}

In the next result, we show that the solution to this quadratic recurrence,
when rescaled in a simple way, gives the number of rooted triangulations
with given genus and number of faces.

\begin{thm}\label{trirec}
The number of rooted triangulations
of genus $g$, with $2n$ faces, is given by
$$F(n,g)=\frac{1}{3n+2}f(n,g),\qquad\qquad (n,g)\in\cS\setminus\{(-1,0),(0,0)\}.$$
\end{thm}

\begin{proof}
Let $F(n,g)$ be the number of rooted triangulations of genus $g$, with $2n$ faces, for $n\ge 1$, $g\ge 0$. These have $3n$ edges, since the sum of the face degrees is twice the number of edges (which also explains why triangulations must have an even number of faces). They also have $n+2-2g$ vertices, which follows from Euler's formula. Thus
\begin{equation}\label{cube}
V_{\{ 3\}}=\sum_{\substack{n\ge 1\\0\le g\le (n+1)/2}} F(n,g) p_3^{2n}w^{n+2-2g}z^{3n},
\end{equation}
where duality has been used to interchange the numbers of vertices and faces.
Now let
\begin{equation*}
\Psi(p_1,p_2,p_3,w,z):=\int V_{\{ 1,2,3\}}(p_1,p_2,p_3,w,z)\frac{dz}{2z}.
\end{equation*}
Then from Theorem~\ref{mapKP} and~(\ref{parRmu}), we know that $\Psi$ satisfies~(\ref{canonKP}), in
which all partials are with respect to $p_1,p_2,p_3$ (so setting $p_i=0$ for $i\ge 4$ is permissible even
before differentiation). Let $\ta$ denote the substitution operator $p_1\mapsto0, p_2\mapsto 0$. Now apply $\ta$ to~(\ref{canonKP}), to obtain
\begin{equation}\label{KPta}
\ta F_{2,2}-\ta F_{3,1}+\frac{1}{12}\ta F_{1^4}+\frac{1}{2}\left(\ta F_{1,1}\right)^2=0,
\end{equation}
where, from Proposition~\ref{cubeconst} and~(\ref{cube}), the corresponding terms for the solution $\Psi$ are given by
\begin{eqnarray*}
\ta\Psi_{2,2}&=&\frac{1}{2}w^2z^2+\!\!\!\!\!\!\!\!
\sum_{\substack{n\ge 1\\0\le g\le (n+1)/2}}\!\!\!\!\!\!\!\!\frac{3n+1}{2}F(n,g) p_3^{2n}w^{n+2-2g}z^{3n+2},\\
\ta\Psi_{3,1}&=&w^2z^2+\!\!\!\!\!\!\!\!
\sum_{\substack{n\ge 1\\0\le g\le (n+1)/2}}\!\!\!\!\!\!\!\!(2n+1)F(n,g) p_3^{2n}w^{n+2-2g}z^{3n+2},\\
\ta\Psi_{1,1}&=&wz+4p_3^2w^2z^4+\!\!\!\!\!\!\!\!
\sum_{\substack{n\ge 1\\0\le g\le (n+1)/2}}\!\!\!\!\!\!\!\!2(3n+2)F(n,g) p_3^{2n+2}w^{n+2-2g}z^{3n+4},\\
\ta\Psi_{1^4}&=&12p_3^2wz^5+ 384p_3^4w^2z^8+\!\!\!\!\!\!\!\!
\sum_{\substack{n\ge 1\\0\le g\le (n+1)/2}}\!\!\!\!\!\!\!\! 8(3n+2)(3n+4)(3n+6)F(n,g) p_3^{2n+4}w^{n+2-2g}z^{3n+8}.
\end{eqnarray*}
Many of the initial terms may be absorbed into the summations above by
using the set $\cS$. For example, 
$$\ta\Psi_{1^4}=
\sum_{(n,g)\in\cS}8(3n+2)(3n+4)(3n+6)F(n,g) p_3^{2n+4}w^{n+2-2g}z^{3n+8}.$$
The result follows by substituting the summation expressions for the four partials into~(\ref{KPta}), equating
the coefficients of $p_3^{2n}w^{n+2-2g}z^{3n+2}$ on both sides of the resulting equation, 
and rescaling to $(3n+2)F(n,g)=f(n,g)$.
\end{proof}

For example, applying~(\ref{trirecence}) recursively, Theorem~\ref{trirec} gives $F(0,0)=1$, and
$$F(1,0)=4,\qquad\qquad F(1,1)=1,\qquad\qquad F(2,0)=32,\qquad\qquad F(2,1)=28,$$
which are consistent with the tables given in~\cite{jv3}. Jason Gao (private communication) has
checked that the recurrence correctly gives the first $40$ terms in genus $0$ and $1$.  The recurrence~(\ref{trirecence}) for triangulations (scaled by $3n+2$ as in Theorem~\ref{trirec}), appears substantially simpler that the one that has appeared in~\cite{g}, but we do not not know of a direct combinatorial argument for this recurrence.

Bender, Gao and Richmond~\cite{bgr} have been able to use this new recurrence to obtain
the explicit asymptotics for triangulations with fixed $n$ and $g$. Previously, the
most explicit asymptotic form for triangulations had been given by Gao~\cite{g},
which involved a scalar $t_g$ depending only on the genus $g$. The asymptotics
for many classes of rooted and unrooted maps in an orientable surface of genus $g$ was
also known up to the same scalar $t_g$. Thus the fact that
the explicit form in~\cite{bgr} determines $t_g$ explicitly means
that  the asymptotics for all these classes of maps are now
known explicitly.

\section*{Acknowledgements}

We thank Jason Gao, Michael Gekhtman, Kevin Purbhoo, Bruce Richmond and Ravi Vakil for helpful suggestions, and Andrei Okounkov for
suggesting that we look at~\cite{mjd} in order to understand the quadratic differential
equations in~\cite{gs}.


\begin{thebibliography}{ELSV}


\bibitem[BGR]{bgr} E.A.Bender, Z. Gao and L.B.Richmond,
\emph{Calculating the asymptotic coefficients $t_g$ using a new recursion for rooted cubic maps},
preprint 2008.

\bibitem[BMS]{bms} M. Bousquet-M\'elou and G. Schaeffer,
\emph{Enumeration of planar constellations},
Adv. \ Appl. \ Math. \ 24 (2000), 337 -- 368. 


\bibitem[G]{g} Zhi-Cheng Gao,
\emph{The number of rooted triangular maps on a surface},
J. Comb. Theory (B) 52 (1991), 236 -- 249.

\bibitem[GJ1]{gj1} I.P. Goulden and D.M. Jackson,
Combinatorial Enumeration,
John Wiley and Sons, New York, 1983 (reprinted by Dover, 2004).

\bibitem[GJ2]{gj2} I.P. Goulden and D.M. Jackson,
\emph{Transitive factorizations into transpositions and holomorphic
mappings on the sphere},
Proc. \ Amer. \ Math. \ Soc.\ 125 (1997), 51 -- 60.

\bibitem[GS]{gs} I.P. Goulden and Luis. G. Serrano,
\emph{A Simple Recurrence for Covers of the Sphere With
Branch Points of Arbitrary Ramification},
Annals \ Comb. \ 10 (2006), 431 -- 441.

\bibitem[GJV]{gjv} I.P. Goulden, D.M. Jackson and R. Vakil,
\emph{Towards the geometry of double Hurwitz numbers},
Advances Math. \ 198 (2005), 43 --92. 

\bibitem[H]{hur} A. Hurwitz,
\emph{Ueber Riemann'sche Fl\"achen mit gegebenen Verzweigungspunkten},
Math. \ Ann. \ 39 (1891), 1 -- 60.

\bibitem[J]{j} D.M. Jackson,
\emph{On an integral representation for the genus series for $2$-cell
embeddings},
Trans. \ Amer. \ Math. \ Soc. \ 344 (1994), 755 -- 772.

\bibitem[JV1]{jv1} D.M. Jackson, T. Visentin,
\emph{A character theoretic approach to embeddings of rooted maps in an
orientable surface of given genus},
Trans. \ Amer. \ Math. \ Soc. \ 322 (1990), 343 -- 363.

\bibitem[JV2]{jv2} D.M. Jackson, T. Visentin,
\emph{Character theory and rooted maps in an
orientable surface of given genus: Face-colored maps},
Trans. \ Amer. \ Math. \ Soc. \ 322 (1990), 365 -- 376.

\bibitem[JV3]{jv3} D.M. Jackson and T.I. Visentin,
An atlas of the smaller maps in orientable and non-orientable surfaces,
CRC Press, 2001.

\bibitem[Ka]{ka} M. Kazarian,
\emph{KP hierarchy for Hodge integrals},
preprint 2007.

\bibitem[KL]{kl} M. Kazarian and S. Lando,
\emph{An algebro-geometric proof of Witten's conjecture},
preprint 2005, math.AG/0601760.

\bibitem[Ko]{ko} M. Kontsevich,
\emph{Intersection theory on the moduli space of curves and
the matrix Airy function},
Comm. \ Math. \ Phys. \ 147 (1992), 1 -- 23.

\bibitem[LZ]{lz} S.K. Lando and A.K. Zvonkin,
Graphs on Surfaces and Their Applications,
Encyclopaedia of Math. Sci., vol 141, Springer-Verlag, Berlin, 2004.


\bibitem[MJD]{mjd} T. Miwa, M. Jimbo and E. Date,
Solitons: Differential Equations, Symmetries and Infinite Dimensional
Algebras, Cambridge University Press, Cambridge, 2000.

\bibitem[Ok]{o} A. Okounkov,
\emph{Toda equations for Hurwitz numbers},
Math. \ Res. \ Letters \ 7 (2000), 447 --453.

\bibitem[Or]{or} A. Orlov,
\emph{Hypergeometric functions as infinite-soliton Tau functions},
Theoretical \ Math. \ Phys.,  146 (2006), 183 --206.

\bibitem[OS]{os} A. Orlov and  D.M. Shcherbin,
\emph{Hypergeometric solutions of soliton equations},
Theoretical \ Math. \ Phys.,  128 (2001), 906 -- 926.

\bibitem[P]{p} R. Pandharipande,
\emph{The Toda equations and the Gromov-Witten theory of the Riemann sphere},
Lett. \ Math. \ Phys. \ 53 (2000), 59 -- 74.

\bibitem[S]{s} R.P. Stanley,
Enumerative Combinatorics, Volume 2,
Cambridge University Press, 1999.

\bibitem[T]{t} W.T. Tutte,
Graph Theory,
Encyclopedia of Math. \ and Applns. \ 21, Addison-Wesley, London, 1984.

\bibitem[W]{w} E. Witten,
\emph{Two dimensional gravity and intersection theory on moduli space},
Survey \ Diff. \ Geom. \ 1 (1991), 243 -- 310.

\end{thebibliography}
\end{document}